\documentclass{aomart}
\usepackage{amssymb}
\numberwithin{equation}{section}
\newtheorem{theorem}{Theorem}[section]
\newtheorem{lemma}[theorem]{Lemma}
\theoremstyle{remark}
\newtheorem*{remark}{Remark}
\newcommand{\norm}[1]{\lVert #1\rVert}
\newcommand{\R}{\mathbb R}
\newcommand{\dd}{\,\mathrm d}
\newcommand{\im}{\operatorname{Im}}
\newcommand{\diver}{\operatorname{div}}
\newcommand{\re}{\operatorname{Re}}
\makeatletter
\let\annals@settitle\@settitle
\renewcommand{\@settitle}{%
  \begingroup
  \def\LARGE{\fontsize{14.4}{18}\selectfont}
  \annals@settitle
  \endgroup}
\AtBeginDocument{\def\@aom@by{}}
\makeatother
\title[Radial blow-up for the Zakharov system]
 {Finite-time blow-up for radial negative-energy solutions of~the Zakharov
 system in two and three dimensions}
\author[Liu]{Jie Liu}
\address{New York University Abu Dhabi, Abu Dhabi, United Arab Emirates}
\email{jl15817@nyu.edu}
\copyrightnote{}
\doinumber{}

\begin{document}
\begin{abstract}
We prove that radial solutions of the Zakharov system in two and three
space dimensions with negative energy and initial data in
$H^2\times H^1\times H^1$ blow up in finite time in both time directions.
This rules out infinite-time grow-up in Merle's dichotomy. The proof
uses a localized virial estimate that controls both the virial and
its dependence on the localization radius, allowing us to propagate a
positive lower bound along shrinking radii and obtain a contradiction
in finite time.
\end{abstract}
\maketitle
\thispagestyle{plain}

\section{Introduction}

The Zakharov system \cite{Zakharov} couples high-frequency Langmuir
oscillations to low-frequency ion-density fluctuations. In its
Hamiltonian formulation on $\R^d$, $d\in\{2,3\}$, the equations are
\begin{equation}\label{eq:zakharov}
 \begin{cases}
 i\partial_t\psi+\Delta\psi=n\psi,\\
 \partial_t n=-\diver v,\\
 \partial_t v=-c_0^2\nabla(n+|\psi|^2).
 \end{cases}
\end{equation}
Here $\psi$ is the complex envelope, $n$ is the real density
fluctuation, $v$ is a real gradient vector field, and $c_0>0$ is the
ion sound speed. All spatial norms and integrals are over $\R^d$.
The mass $M=\norm{\psi}_2^2$ and energy (Hamiltonian)
\begin{equation}\label{eq:energy}
 E=\norm{\nabla\psi}_2^2+\int n|\psi|^2
            +\tfrac12\norm{n}_2^2+\tfrac1{2c_0^2}\norm{v}_2^2
\end{equation}
are conserved. We work in the regular Cauchy class considered in
\cite[p.~434]{Merle}, with solutions continuous in
$H^2\times H^1\times H^1$ on their maximal intervals of existence.
Uniqueness preserves radiality.

The formal limit $c_0\to\infty$ gives $n=-|\psi|^2$ and the
focusing cubic nonlinear Schr\"odinger equation. Schochet and Weinstein
\cite{SW} justified the limit for suitable data on a common local
time interval. For the limiting equation, Glassey's variance argument
\cite{Glassey} gives finite-time blow-up for negative-energy data of
finite variance. Ogawa and Tsutsumi \cite{OTNLS} treat radial data
by a localized virial argument. In the Zakharov system, the
virial also contains the wave momentum. Its effect on the blow-up
argument is described in Section~\ref{sec:virial}.

For $d=2,3$, Merle \cite[Prop.~4.1 and Cor.~4.2]{Merle} proved that
a radial negative-energy solution either blows up in finite time or
exists globally with
\[
 \norm{\psi(t)}_{H^1}+\norm{n(t)}_2+\norm{v(t)}_2
 \longrightarrow\infty\qquad(t\to\infty).
\]
We call the second alternative \emph{grow-up}. Merle
\cite[Sect.~V]{Merle} asked whether this alternative can occur.
His three-dimensional finite-time result
\cite[Prop.~5.1]{Merle} assumes that the time integral of $\norm{\nabla\psi}_2^2$
dominates that of a superlinear power of $\norm{n}_2^2$.
Gan and Guo \cite[Thm.~3.1]{GG10} obtained a finite-time blow-up criterion
for a generalized three-dimensional Zakharov system that includes the
classical system when the added power term vanishes. In the classical
case, they assume negative energy, a positive initial localized virial,
and a time-integrated dominance condition of the same form as Merle's.
Guo, Nakanishi, and Wang \cite[Thm.~1.1]{GNW} obtained a
three-dimensional radial dichotomy below the ground-state
mass--energy threshold. Its blow-up side includes the negative-energy
case and permits both finite-time blow-up and grow-up. Guo and
Nakanishi \cite[Thm.~1.2 and Rem.~1]{GN26} studied the dynamics
near the three-dimensional ground-state threshold; their blow-up
class also allows both possibilities.

Constructions of blow-up profiles address a different question.
In dimension two, Glangetas and Merle \cite{GM1} constructed a
self-similar family. Ghoul, Liu, and Masmoudi \cite[Thm.~1.2]{GhLM}
study the stability of its small-parameter branch. Yin
\cite[Thm.~1.1]{Yin26} constructs particular smooth radial
asymptotically self-similar blow-up solutions in three dimensions.
The present result instead gives a finite-time blow-up criterion
for the entire radial negative-energy class.

\begin{theorem}\label{thm:main}
Let $d\in\{2,3\}$ and $c_0>0$. Suppose
$(\psi_0,n_0,v_0)\in H^2\times H^1\times H^1$ are radial and have
negative Hamiltonian. Then the maximal existence interval of the
regular solution of \eqref{eq:zakharov} with these data is finite
in both time directions.
\end{theorem}

\begin{remark}
The conclusion also holds for radial initial data in
$H^1\times L^2\times L^2$, in the same Hamiltonian formulation.
This follows by smooth approximation and the energy-space local
theory of Bourgain and Colliander \cite{BC}.
Moreover, the continuation criterion and the conserved energy imply
that $\norm{\nabla\psi(t)}_2\to\infty$ at either finite endpoint
of the maximal lifespan.
\end{remark}

\newpage
\subsection{Outline of the proof}\label{sec:outline}

\paragraph{\textbf{Fixed-radius virial.}}
For negative energy, Merle's signed local virial has a strictly positive
time derivative at a sufficiently large fixed radius \cite[Prop.~4.1]{Merle}.
To see the effect of the Zakharov coupling, use the multiplier and weight
\[
 F_R(x)=\frac{x}{\sqrt{1+|x|^2/R^2}},\qquad
 \theta_R=\frac{|F_R|^2}{R^2}.
\]
Write $\partial_r$ for the radial derivative, where $r=|x|$. For a
radial cubic Schr\"odinger solution $u$, the signed momentum
$P_R=-\im\int\bar u F_R\cdot\nabla u$ satisfies
\[
 \frac{|P_R|^2}{R^2}
 \le \norm{u(0)}_2^2\int\theta_R|\partial_r u|^2.
\]
Thus a positive lower bound for $P_R'$ by the weighted kinetic term
yields a quadratic Riccati inequality at fixed radius. For Zakharov,
$Y_R$ in \eqref{Y} also contains $-c_0^{-2}\int n F_R\cdot v$.
Both factors evolve in time; bounding their product by positive
quadratic terms in the virial derivative gives only a linear
differential comparison, compatible with infinite-time grow-up.

\paragraph{\textbf{Control at varying radii.}}
Merle's multiplier $\nabla\phi_R$ equals $x$ on $|x|<R$
\cite[Sect.~IV]{Merle}. There,
$d-\Delta\phi_R=0$ and $\partial_R\nabla\phi_R=0$, while the wave
momentum still contains the integral of $n\,x\cdot v$ there. Thus the
density term in the virial derivative gives no bound for the interior
contribution to $Y_R/R$. Our choice of multiplier instead retains
a weighted density square, which controls the wave momentum and
its radius derivative through the comparisons
\[
 d-\diver F_R\ge \tfrac d2\theta_R,
 \qquad \frac{|F_R|}{R}=\sqrt{\theta_R},
 \qquad |\partial_R F_R|=\theta_R^{3/2}\le\sqrt{\theta_R}.
\]
After absorption of the quartic loss in Section~\ref{sec:proof},
the same positive quadratic terms control both quantities:
for $0<\nu<c_0/2$, Section~\ref{sec:shrinking} gives
\[
 \left(\tfrac{c_0}{2}-\nu\right)\frac{|Y_R|}{R}
     +\nu|\partial_RY_R|
 \le \tfrac12\int\theta_R|\partial_r\psi|^2
       +\tfrac14\int\theta_R n^2
       +\tfrac1{2c_0^2}\norm{v}_2^2+C.
\]
Here $C=C(d,c_0,M)>0$. Along shrinking radii,
\[
 (\partial_t-\nu\partial_R)Y_R
 \ge \frac{c_0/2-\nu}{R}Y_R
      -C\bigl(1+R^{-2}+R^{-2(d-1)}\bigr).
\]

\paragraph{\textbf{The finite-time contradiction.}}
For a suitable choice of $\nu$, multiplication by $R^{4d-5}$ and
subtraction of an explicit function of $R$ give a corrected virial
$\widetilde Y$. The fixed-radius drift makes it positive at some
$(t_0,R_0)$ if the solution lasts until $t_0$, so that
\[
 (\partial_t-\nu\partial_R)\widetilde Y\ge0,\qquad
 \widetilde Y(t_0,R_0)>0,\qquad
 \widetilde Y(t,0)=0.
\]
Along $R=R_0-\nu(t-t_0)$, the positive lower bound persists until
$t=t_0+R_0/\nu$, contradicting the zero boundary value.
Thus the solution cannot continue through this finite time.

\section{The localized virial identity}\label{sec:virial}

For a solution $u$ of the cubic Schr\"odinger equation and a bounded,
smooth, nonnegative weight $\Phi$ with bounded derivatives, the
localized variance satisfies
\[
 \frac{\dd}{\dd t}\int\Phi|u|^2
     =2\im\int\bar u\,\nabla\Phi\cdot\nabla u.
\]
Ogawa and Tsutsumi \cite{OTNLS} use radiality and a suitable localized
weight to obtain a uniformly negative second derivative at a fixed
large radius. Integrating
twice contradicts the nonnegativity of the localized variance.
For Zakharov, the same continuity identity holds for $|\psi|^2$.
In the time derivative of the Schr\"odinger momentum, the term
$-\int|\psi|^2\nabla\Phi\cdot\nabla n$ combines with the wave term
$-\int n\nabla\Phi\cdot\nabla|\psi|^2$ to give
$\int (\Delta\Phi)n|\psi|^2$. The sum of the two momenta is half
the time derivative of
\[
 \int\Phi|\psi(t)|^2+
      \frac{2}{c_0^2}\int_0^t\int n(\tau)\nabla\Phi\cdot v(\tau)
                     \,\dd\tau ,
\]
whose wave correction has no fixed sign. A related fixed-radius argument
of Nguyen and Dinh \cite{ND26} treats a quadratic Schr\"odinger
system. There the conserved charge bounds the undifferentiated
factors in the two momenta, leading to a Riccati inequality.

For $R>0$, define the multiplier and its associated weight by
\begin{equation}\label{eq:multiplier}
 F_R(x)=\frac{x}{\sqrt{1+|x|^2/R^2}},
 \qquad
 \theta_R(x)=\frac{|F_R(x)|^2}{R^2}=\frac{|x|^2}{R^2+|x|^2}.
\end{equation}
The signed localized virial is
\begin{equation}\label{Y}
 Y_R(t)=-\im\int\bar\psi F_R\cdot\nabla\psi
                   -c_0^{-2}\int n F_R\cdot v .
\end{equation}
The Jacobian $DF_R$ is symmetric because $F_R$ is a gradient.
The identity below is Merle's local virial identity
\cite[Eq.~(3.27)]{Merle}, written for the negative momentum
and regrouped using \eqref{eq:energy}. The properties of $F_R$ used
to estimate this identity are established in Section~\ref{sec:proof}.

\begin{lemma}\label{lem:virial}
For a regular solution of \eqref{eq:zakharov}, $(t,R)\mapsto Y_R(t)$
is $C^1$ for $t$ in the lifespan and $R>0$. Writing $I$ for the
identity matrix, we have
\begin{align}
 \partial_tY_R={}&-dE+
   \re\int[(dI-2DF_R)\nabla\psi]\cdot\overline{\nabla\psi}
     +\int(d-\diver F_R)(\tfrac12 n^2+n|\psi|^2)\notag\\
 &+c_0^{-2}\int v\cdot
   [\tfrac12(d+\diver F_R)I-DF_R]v
     +\tfrac12\int\Delta(\diver F_R)|\psi|^2. \label{vir}
\end{align}
\end{lemma}

\begin{proof}
Since $H^2\hookrightarrow L^\infty$ for $d=2,3$, the equations
give continuity of $\partial_t\psi,\partial_t n,\partial_t v$
into $L^2$.
For the Schr\"odinger momentum, set
$\delta_h\psi=(\psi(t+h)-\psi(t))/h$. Integration by parts gives
\[
 \int\overline{\psi(t)}F_R\cdot\nabla\delta_h\psi
 =-\int\bigl(F_R\cdot\nabla\overline{\psi(t)}
        +(\diver F_R)\overline{\psi(t)}\bigr)\delta_h\psi.
\]
Since $\delta_h\psi\to\partial_t\psi$ in $L^2$ and
$\nabla\psi(t+h)\to\nabla\psi(t)$ in $L^2$, passage to the limit
in the time difference quotient gives
\begin{align*}
 \frac{\dd}{\dd t}\im\int\bar\psi F_R\cdot\nabla\psi
 &=2\im\int\overline{\partial_t\psi}\,F_R\cdot\nabla\psi
   -\im\int(\diver F_R)\bar\psi\,\partial_t\psi\\
 &=-2\re\int\Delta\bar\psi F_R\cdot\nabla\psi
   -\re\int(\diver F_R)\bar\psi\,\Delta\psi
   -\int|\psi|^2 F_R\cdot\nabla n .
\end{align*}
The kinetic integrals satisfy
\begin{align*}
 -2\re\int\Delta\bar\psi F_R\cdot\nabla\psi
 &=2\re\int(DF_R\nabla\psi)\cdot\overline{\nabla\psi}
   -\int(\diver F_R)|\nabla\psi|^2,\\
 -\re\int(\diver F_R)\bar\psi\,\Delta\psi
 &=\int(\diver F_R)|\nabla\psi|^2
   -\tfrac12\int\Delta(\diver F_R)|\psi|^2 .
\end{align*}
For the wave momentum, use $\partial_i v_j=\partial_j v_i$ to get
\begin{align*}
 -\int(\diver v) F_R\cdot v
 &=\int v\cdot(DF_R)v +\tfrac12\int F_R\cdot\nabla|v|^2\\
 &=\int v\cdot[DF_R-\tfrac12(\diver F_R)I]v .
\end{align*}
Consequently
\begin{align*}
 c_0^{-2}\frac{\dd}{\dd t}\int n F_R\cdot v
 &=c_0^{-2}\int v\cdot[DF_R-\tfrac12(\diver F_R)I]v\\
 &\quad+\tfrac12\int(\diver F_R)n^2
                 -\int n F_R\cdot\nabla|\psi|^2 .
\end{align*}
The two mixed terms add to $\int(\diver F_R)n|\psi|^2$.
Taking the negative of the sum and using \eqref{eq:energy} gives
\eqref{vir}. All the integrations by parts are justified by
spatial approximation in $H^2\times H^1\times H^1$.
The coefficients are bounded at each positive radius, locally
uniformly in $R$, and pointwise continuous in $R$. Thus the right
side of \eqref{vir} is jointly continuous in $(t,R)$.

Direct differentiation gives
\begin{equation}\label{eq:radius}
\begin{aligned}[b]
 \partial_R F_R&=\frac{|x|^2 x}{(R^2+|x|^2)^{3/2}},\qquad
 |\partial_R F_R|=\theta_R^{3/2}\le\sqrt{\theta_R},\\
 \partial_R^2 F_R&=-\frac{3R}{R^2+|x|^2}\partial_R F_R.
\end{aligned}
\end{equation}
The last formula gives $\norm{\partial_R^2F_R}_\infty\le3/R$.
Thus $F_R$ is $C^1$ as an $L^\infty$-valued function of $R>0$.
In \eqref{Y} one can differentiate $F_R$ under the integrals,
giving the jointly continuous $\partial_R Y_R$.
\end{proof}

\section{Weighted bounds for the localized virial}\label{sec:proof}

The density square and the mixed product in \eqref{vir} have the
same coefficient $d-\diver F_R$. Its lower bound by a multiple of
$\theta_R$ supplies the density control in Section~\ref{sec:shrinking}.
Its upper bound by the same weight allows us to control the quartic
loss by a radial estimate.

For $x\ne0$, put $r=|x|$ and $e_r=x/r$. The radial direction is
$\mathbb R e_r$, and the tangential directions are $e_r^\perp$, the
tangent space of the sphere of radius $r$ at $x$. Differentiating
$F_R=x(1+r^2/R^2)^{-1/2}$ in a direction $\xi\in\R^d$ gives
\begin{align*}
 DF_R(x)\xi
 &=\frac{\xi}{\sqrt{1+r^2/R^2}}
   -\frac{(x\cdot\xi)x}{R^2(1+r^2/R^2)^{3/2}}\\
 &=\sqrt{1-\theta_R}\bigl(\xi-\theta_R(e_r\cdot\xi)e_r\bigr).
\end{align*}
In particular,
\[
 DF_R e_r=(1-\theta_R)^{3/2}e_r,\qquad
 DF_R\xi=\sqrt{1-\theta_R}\,\xi\quad (\xi\perp e_r).
\]
The radial eigenvalue has multiplicity one and the tangential eigenvalue
has multiplicity $d-1$. Taking the trace gives
\[
 d-\diver F_R
 =1-(1-\theta_R)^{3/2}+(d-1)(1-\sqrt{1-\theta_R}).
\]
For $0\le\theta\le1$, the inequality $(1-\theta)^{3/2}\le1-\theta$
and the following identities give
\[
 \frac{\theta}{2}\le\frac{\theta}{1+\sqrt{1-\theta}}
     =1-\sqrt{1-\theta}\le\theta,\qquad
 \theta\le 1-(1-\theta)^{3/2}
     =\frac32\int_0^\theta\sqrt{1-s}\,\dd s\le\frac32\theta .
\]
It follows that
\begin{equation}\label{eq:defect}
 \tfrac d2\theta_R\le d-\diver F_R\le 2d\theta_R.
\end{equation}
The lower bound also holds inside the localization ball. In view of
$|F_R|/R=\sqrt{\theta_R}$ and \eqref{eq:radius}, it supplies the
same density weight for controlling the wave momentum and its radius
derivative.

For radial $\psi$, we have $\nabla\psi=(\partial_r\psi)e_r$.
The Schr\"odinger kinetic coefficient in \eqref{vir} satisfies
\begin{equation}\label{eq:radialkinetic}
 \begin{aligned}
 d-2(1-\theta_R)^{3/2}
 &=d-2+2\bigl(1-(1-\theta_R)^{3/2}\bigr)\\
 &\ge d-2+2\theta_R
 \ge\begin{cases}
 2\theta_R,&d=2,\\
 1,&d=3 .
 \end{cases}
 \end{aligned}
\end{equation}
For radial $v$, the corresponding wave coefficient is
\begin{align*}
 \tfrac12(d+\diver F_R)-(1-\theta_R)^{3/2}
 &=\tfrac d2+\tfrac12\sqrt{1-\theta_R}(d-2+\theta_R)\\
 &\ge\tfrac d2.
\end{align*}
Finally, $F_R=R F_1(\cdot/R)$ gives
$\norm{\Delta\diver F_R}_\infty=R^{-2}\norm{\Delta\diver F_1}_\infty$.

The following estimate bounds the quartic loss with the weight in
\eqref{eq:defect}. Its proof uses the radial integration argument of
\cite[Lemmas~4.3 and~4.4]{Merle}.

\begin{lemma}\label{lem:radial}
For radial $u\in H^1(\R^d)$ and $R>0$,
\begin{equation}\label{quartic}
 \int\theta_R|u|^4\le
 \begin{cases}
 \displaystyle\frac{\norm{u}_2^3}{\pi R}
           \norm{\sqrt{\theta_R}\,\partial_r u}_2,&d=2,\\[5pt]
 \displaystyle\frac{\norm{u}_2^3}{2\pi R^2}
           \norm{\partial_r u}_2,&d=3.
 \end{cases}
\end{equation}
\end{lemma}

\begin{proof}
For smooth compactly supported radial $u$, write $u(r)$ and
$\theta_R(r)$ for the radial profiles. Since $\theta_R'\ge0$,
\begin{align*}
 \theta_R(r)|u(r)|^2
 &=-\int_r^\infty\theta_R'(s)|u(s)|^2\dd s
   -2\re\int_r^\infty\theta_R(s)\bar u(s)u'(s)\dd s\\
 &\le 2\int_r^\infty\theta_R(s)|u(s)u'(s)|\dd s .
\end{align*}
In dimension two, $\sqrt{\theta_R(s)}\le s/R$ and Cauchy--Schwarz
with measure $s\,\dd s$ give
\[
 2\int_r^\infty\theta_R(s)|u(s)u'(s)|\,\dd s
 \le \frac2R\int_0^\infty s|u(s)|\sqrt{\theta_R(s)}|u'(s)|\,\dd s
 \le \frac{\norm{u}_2}{\pi R}\norm{\sqrt{\theta_R}\,\partial_r u}_2.
\]
In dimension three, $\theta_R(s)\le s^2/R^2$ similarly gives
\[
 2\int_r^\infty\theta_R(s)|u(s)u'(s)|\,\dd s
 \le \frac2{R^2}\int_0^\infty s^2|u(s)u'(s)|\,\dd s
 \le \frac{\norm{u}_2}{2\pi R^2}\norm{\partial_r u}_2 .
\]
Multiplying by $|u|^2$ and integrating in space gives \eqref{quartic}
for smooth compactly supported radial $u$. Approximation in $H^1$
and the embedding $H^1\hookrightarrow L^4$ extend the estimate to
all radial $u\in H^1$.
\end{proof}

Reserve a quarter of the density square before estimating the mixed
term. By $\tfrac12 n^2+n|\psi|^2\ge\tfrac14n^2-|\psi|^4$ and
\eqref{eq:defect},
\[
 \int(d-\diver F_R)(\tfrac12n^2+n|\psi|^2)
       \ge \tfrac14\int\theta_R n^2-2d\int\theta_R|\psi|^4 .
\]
By \eqref{quartic} and Young's inequality,
\begin{align*}
 2d\int\theta_R|\psi|^4
 &\le \frac{4M^{3/2}}{\pi R}\norm{\sqrt{\theta_R}\partial_r\psi}_2
 \le \int\theta_R|\partial_r\psi|^2+\frac{4M^3}{\pi^2R^2}
       &&(d=2),\\
 2d\int\theta_R|\psi|^4
 &\le \frac{3M^{3/2}}{\pi R^2}\norm{\partial_r\psi}_2
 \le \tfrac12\norm{\partial_r\psi}_2^2+\frac{9M^3}{2\pi^2R^4}
       &&(d=3).
\end{align*}
Combining these estimates with \eqref{vir} and
\eqref{eq:radialkinetic}, for $E<0$ we obtain
\begin{align}
 \partial_tY_R\ge {}&d|E|
   +\tfrac12\int\theta_R|\partial_r\psi|^2
   +\tfrac14\int\theta_R n^2+\tfrac1{2c_0^2}\norm{v}_2^2
   -C\bigl(R^{-2}+R^{-2(d-1)}\bigr). \label{fixed}
\end{align}
Here and below $C>0$ can be chosen once, depending only on $d$,
$c_0$, and $M$, large enough for both \eqref{fixed} and
\eqref{eq:comparison} below.

\section{Exclusion of grow-up}\label{sec:shrinking}

The positive terms in \eqref{fixed} control both $Y_R/R$ and
$\partial_RY_R$. By \eqref{eq:multiplier}, \eqref{eq:radius}, and
Cauchy--Schwarz, each of $|Y_R|/R$ and $|\partial_R Y_R|$ is bounded by
\begin{equation}\label{dual}
 \sqrt M\,\norm{\sqrt{\theta_R}\partial_r\psi}_2
    +c_0^{-2}\norm{\sqrt{\theta_R}n}_2\norm{v}_2 .
\end{equation}
For $0<\nu<c_0/2$, multiply the two bounds by $c_0/2-\nu$
and $\nu$, respectively. Adding and applying Young's inequality gives
\begin{align*}
 \left(\tfrac{c_0}{2}-\nu\right)\frac{|Y_R|}{R}
       +\nu|\partial_R Y_R|
 &\le \tfrac12\int\theta_R|\partial_r\psi|^2
       +\tfrac14\int\theta_R n^2
       +\tfrac1{2c_0^2}\norm{v}_2^2+\tfrac{c_0^2M}{8}.
\end{align*}
Substitution into \eqref{fixed} yields the differential inequality
along a radius moving inward at speed $\nu$:
\begin{equation}\label{eq:comparison}
 (\partial_t-\nu\partial_R)Y_R
 \ge d|E|+\frac{c_0/2-\nu}{R}Y_R
    -C\bigl(1+R^{-2}+R^{-2(d-1)}\bigr).
\end{equation}

\begin{proof}[Proof of Theorem~\ref{thm:main}]
Set $\nu=c_0/(8(d-1))$ and let $k=4d-5$, so that
$c_0/2-\nu=k\nu$. Multiplying \eqref{eq:comparison} by $R^k$ gives
\begin{align}
 (\partial_t-\nu\partial_R)(R^kY_R)
 &=R^k\left((\partial_t-\nu\partial_R)Y_R-\frac{k\nu}{R}Y_R\right)\notag\\
 &\ge -C\bigl(R^k+R^{k-2}+R^{k-2(d-1)}\bigr). \label{barrier}
\end{align}
All three powers on the right are integrable at zero for $d=2,3$.
With the same constant $C$, define the corrected virial for $R>0$ by
\[
 \widetilde Y(t,R)=R^kY_R(t)-\frac C\nu\int_0^R
          \bigl(r^k+r^{k-2}+r^{k-2(d-1)}\bigr)\dd r .
\]
Differentiating along the shrinking radius and using
\eqref{barrier}, we obtain
\begin{align}
 (\partial_t-\nu\partial_R)\widetilde Y
 &=(\partial_t-\nu\partial_R)(R^kY_R)
   +C\bigl(R^k+R^{k-2}+R^{k-2(d-1)}\bigr)\notag\\
 &\ge0.\label{eq:transport}
\end{align}
Moreover, the bound $|F_R|\le R$ gives
\[
 |R^kY_R(t)|\le R^{k+1}
 \left(\sqrt M\,\norm{\nabla\psi(t)}_2+
 c_0^{-2}\norm{n(t)}_2\norm{v(t)}_2\right).
\]
The factor in parentheses is bounded on each compact subinterval of
the lifespan. The integral in $\widetilde Y$ also tends to zero as
$R\downarrow0$. Thus $\widetilde Y$ extends jointly continuously to
the zero-radius boundary, with
$\widetilde Y(t,0)=0$.

Let $T_+$ be the forward endpoint of the maximal lifespan. Choose
$R_0\ge1$ so large that $C(R_0^{-2}+R_0^{-2(d-1)})\le |E|$.
At that fixed radius, \eqref{fixed} gives
\begin{equation}\label{eq:drift}
 Y_{R_0}(t)\ge Y_{R_0}(0)+(d-1)|E|t\qquad (0\le t<T_+).
\end{equation}
Choose $t_0\ge0$ such that
\[
 R_0^k\bigl(Y_{R_0}(0)+(d-1)|E|t_0\bigr)
 >\frac C\nu\int_0^{R_0}\bigl(r^k+r^{k-2}+r^{k-2(d-1)}\bigr)\dd r .
\]
If $T_+\le t_0$, then $T_+$ is finite. Otherwise,
$\widetilde Y(t_0,R_0)>0$. Put $T_*=t_0+R_0/\nu$. Integrating
\eqref{eq:transport} along $R=R_0-\nu(t-t_0)$ gives
\begin{equation}\label{eq:threshold}
 \widetilde Y(t,R_0-\nu(t-t_0))\ge\widetilde Y(t_0,R_0)>0,
 \qquad t_0\le t<\min\{T_+,T_*\}.
\end{equation}
Suppose $T_+>T_*$. Then $T_*$ lies strictly inside the lifespan,
so joint continuity at $(T_*,0)$ and \eqref{eq:threshold} give
\[
 0=\widetilde Y(T_*,0)
   =\lim_{t\uparrow T_*}\widetilde Y(t,R_0-\nu(t-t_0))
   \ge\widetilde Y(t_0,R_0)>0,
\]
a contradiction. Hence $T_+\le T_*$.
Finally,
\[
 (\psi(t),n(t),v(t))\longmapsto
       (\overline{\psi(-t)},n(-t),-v(-t))
\]
preserves \eqref{eq:zakharov}, radiality, the Hamiltonian, and the
regular Cauchy class. The forward bound applied to the time-reversed
solution gives the backward assertion.
\end{proof}

\begin{samepage}
\section*{AI use declaration}
This work grew out of the author's joint research on self-similar
Zakharov blow-up \cite{GhLM} and his unpublished notes on
negative-energy blow-up near the two-dimensional critical mass.
The author identified the obstruction caused by the wave momentum
and formulated the characteristic interpretation of the proof:
propagation of a positive lower bound to the zero-radius boundary.
Under the author's direction, OpenAI Codex assisted in developing
the varying-radius estimate and integrating-factor correction in
Section~\ref{sec:shrinking}, and in preparing and editing the
manuscript. The author takes full responsibility for the final
manuscript.\par
\end{samepage}

\bibliographystyle{aomplain}
\bibliography{references}

@misc{GhLM,
  author = {Ghoul, Tej-eddine and Liu, Jie and Masmoudi, Nader},
  title = {Asymptotic stability of self-similar blow-up for the two-dimensional {Zakharov} system},
  year = {2026},
  howpublished = {Preprint},
  arxiv = {2611.16731}
}

@article{Glassey,
  author = {Glassey, Robert T.},
  title = {On the blowing up of solutions to the {Cauchy} problem for nonlinear {Schr\"odinger} equations},
  journal = {J. Math. Phys.},
  year = {1977},
  volume = {18},
  number = {9},
  pages = {1794--1797},
  doi = {10.1063/1.523491}
}

@article{GM1,
  author = {Glangetas, Laurent and Merle, Frank},
  title = {Existence of self-similar blow-up solutions for {Zakharov} equation in dimension two. {I}},
  journal = {Comm. Math. Phys.},
  year = {1994},
  volume = {160},
  number = {1},
  pages = {173--215},
  doi = {10.1007/BF02099792},
  mrnumber = {1262194}
}

@misc{GN26,
  author = {Guo, Zihua and Nakanishi, Kenji},
  title = {Global dynamics above the ground state energy for the {3D} {Zakharov} system},
  year = {2026},
  howpublished = {Preprint},
  arxiv = {2608.14993}
}

@article{GNW,
  author = {Guo, Zihua and Nakanishi, Kenji and Wang, Shuxia},
  title = {Global dynamics below the ground state energy for the {Zakharov} system in the {3D} radial case},
  journal = {Adv. Math.},
  year = {2013},
  volume = {238},
  pages = {412--441},
  doi = {10.1016/j.aim.2013.02.008}
}

@article{Merle,
  author = {Merle, Frank},
  title = {Blow-up results of viriel type for {Zakharov} equations},
  journal = {Comm. Math. Phys.},
  year = {1996},
  volume = {175},
  number = {2},
  pages = {433--455},
  doi = {10.1007/BF02102415}
}

@article{SW,
  author = {Schochet, Steven H. and Weinstein, Michael I.},
  title = {The nonlinear {Schr\"odinger} limit of the {Zakharov} equations governing {Langmuir} turbulence},
  journal = {Comm. Math. Phys.},
  year = {1986},
  volume = {106},
  number = {4},
  pages = {569--580},
  doi = {10.1007/BF01463396}
}

@article{Zakharov,
  author = {Zakharov, V. E.},
  title = {Collapse of {Langmuir} waves},
  journal = {Soviet Phys. JETP},
  year = {1972},
  volume = {35},
  number = {5},
  pages = {908--914},
  url = {https://www.jetp.ras.ru/cgi-bin/e/index/e/35/5/p908?a=list}
}

@article{OTNLS,
  author = {Ogawa, Takayoshi and Tsutsumi, Yoshio},
  title = {Blow-up of {$H^1$} solution for the nonlinear {Schr\"odinger} equation},
  journal = {J. Differential Equations},
  year = {1991},
  volume = {92},
  number = {2},
  pages = {317--330},
  doi = {10.1016/0022-0396(91)90052-B}
}

@misc{Yin26,
  author = {Yin, Yuchen},
  title = {Smooth asymptotically self-similar blow up for the scalar {Zakharov} system in dimension {$3$}},
  year = {2026},
  howpublished = {Preprint},
  arxiv = {2611.12052}
}

@misc{ND26,
  author = {Nguyen, Ngoc Uyen Cong and Dinh, Van Duong},
  title = {Finite-time blow-up for the four-dimensional mass-critical quadratic nonlinear {Schr\"odinger} system without mass resonance},
  year = {2026},
  howpublished = {Preprint},
  arxiv = {2608.18453}
}

@article{GG10,
  author = {Gan, Zaihui and Guo, Boling},
  title = {Finite-time blowup for {Zakharov} system with combined power-type nonlinearities},
  journal = {Appl. Anal.},
  year = {2010},
  volume = {89},
  number = {5},
  pages = {645--661},
  doi = {10.1080/00036811003627476}
}

@article{BC,
  author = {Bourgain, Jean and Colliander, James},
  title = {On wellposedness of the {Zakharov} system},
  journal = {Internat. Math. Res. Notices},
  year = {1996},
  number = {11},
  pages = {515--546},
  doi = {10.1155/S1073792896000359}
}
\end{document}